# SHRINKAGE PRIORS FOR BAYESIAN PREDICTION

By Fumiyasu Komaki

*University of Tokyo*

We investigate shrinkage priors for constructing Bayesian predictive distributions. It is shown that there exist shrinkage predictive distributions asymptotically dominating Bayesian predictive distributions based on the Jeffreys prior or other vague priors if the model manifold satisfies some differential geometric conditions. Kullback–Leibler divergence from the true distribution to a predictive distribution is adopted as a loss function. Conformal transformations of model manifolds corresponding to vague priors are introduced. We show several examples where shrinkage predictive distributions dominate Bayesian predictive distributions based on vague priors.

**1. Introduction.** Suppose that we have a set of independent observations $x^{(N)} = (x(1), x(2), \ldots, x(N))$ from a distribution with density $p(x|\theta)$ that belongs to a model $\{p(x|\theta)|\theta = (\theta^1, \theta^2, \ldots, \theta^d) \in \Theta\}$. An unobserved variable $y := x(N+1)$ from the same distribution $p(y|\theta)$ is predicted by using a predictive density $\hat{p}(y; x^{(N)})$.

We adopt the Kullback–Leibler divergence $D\{p(y|\theta), \hat{p}(y; x^{(N)})\} := \int p(y|\theta) \log\{p(y|\theta)/\hat{p}(y; x^{(N)})\} \, dy$, which has a natural information theoretic meaning, as a loss function. We evaluate the performance of predictive distributions by using the risk function $\mathrm{E}[D(p, \hat{p})|\theta] = \int p(x^{(N)}|\theta) \int p(y|\theta) \log\{p(y|\theta)/\hat{p}(y; x^{(N)})\} \, dy \, dx^{(N)}$.

A widely used method to construct a predictive density is to use a plug-in density $p(y|\hat{\theta}(x^{(N)}))$, where $\hat{\theta}(x^{(N)})$ is an appropriate estimator of $\theta$. However, Bayesian predictive densities

$$p_\pi(y|x^{(N)}) = \frac{\int p(y|\theta) p(x^{(N)}|\theta) \pi(\theta) \, d\theta}{\int p(x^{(N)}|\bar{\theta}) \pi(\bar{\theta}) \, d\bar{\theta}}$$

have better performance than plug-in distributions in many examples [1, 12, 18].

In the present paper we investigate the use of shrinkage priors for constructing Bayesian predictive distributions asymptotically dominating those based on improper vague priors such as the Jeffreys prior.

There exist many studies on shrinkage estimators where elliptic operators including the Laplacian and (super) harmonic functions play important roles [5, 7, 10, 14, 23, 24, 25].

Recently, several results suggesting that shrinkage priors are useful for various prediction problems have been obtained; see [13, 19] for the normal model and [21] for the Poisson model. Further studies on more general models are required.

In the present paper construction methods for shrinkage priors are introduced and properties of them are investigated from the viewpoint of information geometry by using the results of previous studies on asymptotic properties of predictive distributions [8, 15, 18, 28]. The model $\{p(x|\theta)|\theta \in \Theta\}$ is regarded as a manifold, and the relation between differential geometric properties of the model manifold and the existence of shrinkage priors is studied. It is shown that there exist useful shrinkage priors if the model manifold satisfies some differential geometric conditions. The geometrical approach is useful to investigate Bayesian methods because prior distributions are naturally regarded as volume elements on model manifolds.

In Section 2 we show that there exist shrinkage predictive distributions asymptotically dominating the Bayesian predictive distribution based on the Jeffreys prior if the model manifold endowed with the Fisher metric satisfies some differential geometric conditions. In Section 3 we introduce conformal transformations of model manifolds corresponding to prior distributions and show that there exist shrinkage predictive distributions asymptotically dominating Bayesian predictive distributions based on various priors if the transformed model manifolds satisfy some differential geometric conditions. In Section 4 we show several examples where shrinkage predictive distributions constructed by using the introduced methods asymptotically or exactly dominate Bayesian predictive distributions based on vague priors.

**2. Shrinkage priors asymptotically dominating the Jeffreys prior.** First, we present some differential geometric notions and notation to be used. In the following, we assume that the model manifold $M$ is a $d(\geq 2)$-dimensional connected and orientable $C^\infty$ manifold. The parameter space $\Theta$ is regarded as a coordinate system of $M$. We use Einstein's summation convention: if an index occurs twice in any one term, once as an upper and once as a lower index, summation over that index is implied.

The Fisher metric tensor is defined by $g_{ij}(\theta) := \mathrm{E}[\partial_i \log p(x|\theta) \partial_j \log p(x|\theta)|\theta]$, where $\partial_i := \partial/\partial\theta^i$. The coefficients of the $\alpha$-connection are defined by $\overset{\alpha}{\Gamma}{}^k_{ij}(\theta) := \Gamma^k_{ij}(\theta) - \frac{\alpha}{2}T_{ijl}(\theta)g^{kl}(\theta)$, where $\Gamma^k_{ij} := \frac{1}{2}(\partial_i g_{jl}(\theta) + \partial_j g_{li}(\theta) -$

$\partial_l g_{ij}(\theta))g^{kl}(\theta)$ are the coefficients of the Riemannian connection, $T_{ijk}(\theta) := \mathrm{E}[\partial_i \log p(x|\theta)\partial_j \log p(x|\theta)\partial_k \log p(x|\theta)|\theta]$ is the skewness tensor, and $g^{ij}$ denotes the $(i,j)$-component of the inverse matrix of $(g_{ij})$; see [2] for details. The $-1$-connection and $1$-connection are called the m-connection and e-connection and their coefficients are denoted by $\overset{m}{\Gamma}{}^k_{ij}$ and $\overset{e}{\Gamma}{}^k_{ij}$, respectively. The $\alpha$-covariant derivative of a vector field $v$ is defined by $\overset{\alpha}{\nabla}_i v^j = \partial_i v^j + \overset{\alpha}{\Gamma}{}^j_{ik} v^k$, and $\overset{0}{\nabla}$ and $\overset{1}{\nabla}$ are denoted by $\nabla$ and $\overset{e}{\nabla}$, respectively.

The Laplacian $\Delta$ on a manifold $(M,g)$ endowed with a Riemannian metric $g_{ij}$ is defined by

$$\Delta f = |g|^{-1/2}\partial_i(|g|^{1/2}g^{ij}\partial_j f) = \nabla_i(g^{ij}\partial_j f),$$

where $f$ is a real function on $M$, and $|g|$ is the determinant of the matrix $(g_{ij})$.

A continuous function $G(\xi,\theta)$ of $\xi$ and $\theta$ on $M \times M - \{\theta = \xi\}$ is called a Green function if it satisfies the following conditions (see, e.g., [4, 27]):

1. $\Delta_\theta G(\xi,\theta) = 0$ for all $\xi \in M$ and $\theta \neq \xi$, where $\Delta_\theta$ denotes the Laplacian with respect to $\theta$.
2. $G(\xi,\theta) \geq 0$.
3. In a neighborhood of $\xi$, $G(\xi,\theta)$ has the singularity

$$G(\xi,\theta) \sim \begin{cases} \{(d-2)\omega_{d-1}\}^{-1} \mathrm{dis}(\xi,\theta)^{-(d-2)}, & d \geq 3, \\ -(1/2\pi)\log \mathrm{dis}(\xi,\theta), & d = 2, \end{cases}$$

where $\omega_{d-1} := 2\pi^{d/2}/\Gamma(d/2)$ is the area of the $(d-1)$-dimensional unit sphere and $\mathrm{dis}(\xi,\theta)$ denotes the Riemannian distance between $\xi$ and $\theta$.
4. There exists a positive number $\delta > 0$ such that $G(\xi,\theta)$ is a bounded function of $\theta \in \{\theta | \theta \in M, \mathrm{dis}(\xi,\theta) > \delta\}$.

When a Green function $G(\xi,\theta)$ exists on $(M,g)$, it is represented by

$$(1) \qquad G(\xi,\theta) = \int_0^\infty \gamma_t(\xi,\theta)\,dt,$$

where $\gamma_t(\xi,\theta)$ is the minimal positive fundamental solution of the heat equation $\partial u(t,\theta)/\partial t = \Delta u(t,\theta)$; see [16, 17].

In the following, we introduce shrinkage priors and evaluate the risk of Bayesian predictive distributions based on the shrinkage priors by using the results of previous studies on asymptotic properties of predictive distributions. Asymptotic expansions of predictive distributions are studied by Vidoni [28], Komaki [18] and Hartigan [15].

THEOREM 1. *Let $(M,g)$ be a model manifold endowed with the Fisher metric. If a Green function $G(\xi,\theta)$ on $(M,g)$ exists, there exist Bayesian* *predictive distributions asymptotically dominating the Bayesian predictive distribution based on the Jeffreys prior $\pi_{\mathrm{J}}(\theta) \propto |g(\theta)|^{1/2}$. In particular, the Bayesian predictive distribution based on the prior $\pi_{\mathrm{G}}(\theta)\, d\theta := G(\xi, \theta)\pi_{\mathrm{J}}(\theta)\, d\theta$, where $\xi \in M$ is an arbitrary fixed point, asymptotically dominates the Bayesian predictive distribution based on the Jeffreys prior.*

PROOF. In the following we assume that $d \geq 2$ since a Green function does not exist on the manifold when $d = 1$; see [4].

First, we assume that the true parameter value $\theta$ is different from $\xi$. The Bayesian predictive density based on a prior $f(\theta)$ can be expanded as

$$
\begin{aligned}
p_f(y|x^{(N)}) = {} & p(y|\hat{\theta}_{\mathrm{mle}}(x^{(N)})) \\
& + \frac{1}{2N} g^{ij}(\hat{\theta}_{\mathrm{mle}})(\partial_i \partial_j p(y|\hat{\theta}_{\mathrm{mle}}) - \overset{\mathrm{m}}{\Gamma}{}^k_{ij} \partial_k p(y|\hat{\theta}_{\mathrm{mle}})) \\
& + \frac{1}{N}\left\{\partial_i \log \frac{f}{\pi_{\mathrm{J}}}(\hat{\theta}_{\mathrm{mle}}) + \frac{1}{2} T_i(\hat{\theta}_{\mathrm{mle}})\right\} g^{ik}(\hat{\theta}_{\mathrm{mle}}) \partial_k p(y|\hat{\theta}_{\mathrm{mle}}) \\
& + o_p(N^{-1}),
\end{aligned}
\tag{2}
$$

where $\hat{\theta}_{\mathrm{mle}}(x^{(N)})$ is the maximum likelihood estimate, $T_i := T_{ijk} g^{jk}$ and the relation $\partial_i \log \pi_{\mathrm{J}} = \partial_i \log |g_{ij}|^{1/2} = \Gamma^j_{ij} = \overset{\mathrm{e}}{\Gamma}{}^j_{ij} + (1/2)T_i$ is used; see [8, 15, 18].

The risk of the Bayesian predictive density $p_f(y|x^{(N)})$ is given by

$$
\begin{aligned}
& \mathrm{E}[D(p(y|\theta), p_f(y|x^{(N)}))|\theta] \\
& = \frac{d}{2N} + \frac{1}{2N^2} g^{ij}\left(\partial_i \log \frac{f}{\pi_{\mathrm{J}}} + \frac{1}{2}T_i\right)\left(\partial_j \log \frac{f}{\pi_{\mathrm{J}}} + \frac{1}{2}T_j\right) \\
& \quad + \frac{1}{N^2} \overset{\mathrm{e}}{\nabla}_i \left\{g^{ij}\left(\partial_j \log \frac{f}{\pi_{\mathrm{J}}} + \frac{1}{2}T_j\right)\right\} \\
& \quad + \text{the terms independent of } f + o(N^{-2});
\end{aligned}
\tag{3}
$$

see [15, 18].

Thus, the difference between the risk of $p_{\pi_{\mathrm{J}}}(y|x^{(N)})$ based on the Jeffrey prior $\pi_{\mathrm{J}}(\theta)$ and the risk of $p_f(y|x^{(N)})$ based on $f(\theta)$ is given by

$$
\begin{aligned}
& \mathrm{E}[D(p(y|\theta), p_{\pi_{\mathrm{J}}}(y|x^{(N)}))|\theta] - \mathrm{E}[D(p(y|\theta), p_f(y|x^{(N)}))|\theta] \\
& = \frac{1}{8N^2} g^{ij} T_i T_j + \frac{1}{2N^2} \overset{\mathrm{e}}{\nabla}_i (g^{ij} T_j) \\
& \quad - \frac{1}{2N^2} g^{ij}\left(\partial_i \log \frac{f}{\pi_{\mathrm{J}}} + \frac{1}{2}T_i\right)\left(\partial_j \log \frac{f}{\pi_{\mathrm{J}}} + \frac{1}{2}T_j\right) \\
& \quad - \frac{1}{N^2} \overset{\mathrm{e}}{\nabla}_i \left\{g^{ij}\left(\partial_j \log \frac{f}{\pi_{\mathrm{J}}} + \frac{1}{2}T_j\right)\right\} + o(N^{-2})
\end{aligned}
\tag{4}
$$

$$= \frac{1}{2N^2} g^{ij} \partial_i \log \frac{f}{\pi_J} \partial_j \log \frac{f}{\pi_J} - \frac{1}{N^2} \frac{\pi_J}{f} \Delta \frac{f}{\pi_J} + o(N^{-2}).$$

Therefore, if $\Delta(f/\pi_J) \leq 0$ and $\partial_i \log(f/\pi_J) \neq 0$, $p_f(y|x^{(N)})$ asymptotically dominates $p_{\pi_J}(y|x^{(N)})$. The prior density $\pi_G(\theta)$ satisfies these conditions.

Next, we assume that $\theta = \xi$.

Since $\int (\bar{\theta}^i - \theta^i)(\bar{\theta}^j - \theta^j) p_f(\bar{\theta}|x^{(N)}) d\bar{\theta} = O_p(N^{-1})$, $\int (\bar{\theta}^i - \theta^i)(\bar{\theta}^j - \theta^j) \times (\bar{\theta}^k - \theta^k) p_f(\bar{\theta}|x^{(N)}) d\bar{\theta} = O_p(N^{-3/2})$ and $\hat{\theta}^i_f - \theta^i = O_p(N^{-1/2})$, where $\hat{\theta}_f(x^{(N)}) = \int \bar{\theta} p_f(\bar{\theta}|x^{(N)}) d\bar{\theta}$, we have

$$\mathrm{E}[D(p(y|\theta), p_f(y|x^{(N)}))|\theta]$$

$$= \int p(x^{(N)}|\theta) \int p(y|\theta) \log \frac{p(y|\theta)}{p_f(y|x^{(N)})} \, dy \, dx^{(N)}$$

$$= -\int p(x^{(N)}|\theta) \int p(y|\theta) \log \left[ \int \left\{ 1 + \frac{\partial_i p(y|\theta)}{p(y|\theta)} (\bar{\theta}^i - \theta^i) \right. \right.$$

$$+ \frac{1}{2} \frac{\partial_i \partial_j p(y|\theta)}{p(y|\theta)} (\bar{\theta}^i - \theta^i)$$

(5)

$$\times (\bar{\theta}^j - \theta^j) \Big\}$$

$$\times p_f(\bar{\theta}|x^{(N)}) d\bar{\theta} \Bigg] dy \, dx^{(N)}$$

$$+ O(N^{-3/2})$$

$$= \frac{1}{2} g_{ij}(\theta) \int p(x^{(N)}|\theta)(\hat{\theta}^i_f(x^{(N)}) - \theta^i)(\hat{\theta}^j_f(x^{(N)}) - \theta^j) \, dx^{(N)} + O(N^{-3/2}).$$

This relation holds for all the three cases $f = \pi_J$, $f = \pi_G$ ($\xi = \theta$) and $f = \pi_G$ ($\xi \neq \theta$); see [22] for details.

When $f = \pi_J$, the asymptotic distribution of $N g_{ij}(\theta)(\hat{\theta}^i_f(x^{(N)}) - \theta^i) \times (\hat{\theta}^j_f(x^{(N)}) - \theta^j)$ is the chi-square distribution with $d$ degrees of freedom since $p_{\pi_J}(\mu|x^{(N)}) \propto (2\pi)^{-d/2} |g_{ij}(\theta)|^{1/2} \exp\{(1/2) g_{ij}(\theta)(\mu^i - \eta^i)(\mu^j - \eta^j)\}(1 + O_p(N^{-1/2}))$, where $\mu^i := \sqrt{N}(\bar{\theta}^i - \theta^i)$, $\eta^i := \sqrt{N}(\hat{\theta}^i - \theta^i)$ and $\hat{\theta}$ is the maximum likelihood estimator based on the observation $x^{(N)}$. Thus, the risk (5) is $d/(2N) + o(N^{-1})$, coinciding with (3).

When $d \geq 3$, the risk (5) with $f = \pi_G$ ($\xi = \theta$) is smaller than that with $f = \pi_J$ on the order of $1/N$ since $p_{\pi_G}(\mu|x^{(N)}) \propto (g_{ij}(\theta)\mu^i \mu^j)^{-(d-2)/2} (2\pi)^{-d/2} \times |g_{ij}(\theta)|^{1/2} \exp\{(1/2) g_{ij}(\theta)(\mu^i - \eta^i)(\mu^j - \eta^j)\}(1 + O_p(N^{-1/2}))$.

When $d = 2$, we can verify that the risk (5) with $f = \pi_G$ ($\xi = \theta$) is smaller than that with with $f = \pi_J$ on the order of $1/(N \log N)$ since $p_{\pi_G}(\mu|x^{(N)}) \propto \{1 - (1/\log N) \log(g_{ij}(\theta)\mu^i \mu^j)\} \times (2\pi)^{-d/2} |g_{ij}(\theta)|^{1/2} \exp\{(1/2) g_{ij}(\theta)(\mu^i - \eta^i)(\mu^j - \eta^j)\}(1 + O_p(\log N)^{-2})$.

Therefore, the Bayesian predictive distribution based on $\pi_{\mathrm{G}}$ asymptotically dominates that based on the Jeffreys prior $\pi_{\mathrm{J}}$. □

From (4) and the relation $(1/2)g^{ij}\partial_i\log(f/\pi_{\mathrm{J}})\partial_j\log(f/\pi_{\mathrm{J}}) - (\pi_{\mathrm{J}}/f)\Delta(f/\pi_{\mathrm{J}}) = -2(\pi_{\mathrm{J}}/f)^{1/2}\Delta(f/\pi_{\mathrm{J}})^{1/2}$, we have the following theorem. For the definition of superharmonic functions on Riemannian manifolds, see, for example, [16]. A $C^2$ function is superharmonic if and only if $\Delta f \leq 0$.

THEOREM 2. *Let $f(\theta)$ be a smooth prior density on a model manifold $(M,g)$ endowed with the Fisher metric. The Bayesian predictive distribution based on $f(\theta)$ asymptotically dominates the Bayesian predictive distribution based on the Jeffreys prior $\pi_{\mathrm{J}}(\theta)$ if and only if $(f/\pi_{\mathrm{J}})^{1/2}$ is a nonconstant positive superharmonic function on $(M,g)$.*

It is known that there exists a Green function associated with the Laplacian $\Delta$ if and only if there exists at least one nonconstant positive superharmonic function [16]. Therefore, the existence of a Green function on $(M,g)$ is necessary (and sufficient) for the existence of positive superharmonic functions on $(M,g)$.

It is proved by Aomoto [4] that there exists a Green function if a complete and simply connected manifold has strictly negative curvature ($d = 2$) or has negative curvature ($d \geq 3$).

The sectional curvature of two-dimensional subspace of the tangent space $T_p$ of $M$ at $p$ spanned by $X$ and $Y$ is defined by $K(X,Y) := (R_{ijkl}X^iY^jY^kX^l)/\{(g_{ik}g_{jl} - g_{jk}g_{il})X^iY^jX^kY^l\}$, where $R_{ijkl}$ is the curvature tensor defined by $R_{ijkl} := (\partial_i\Gamma^m_{jk} - \partial_j\Gamma^m_{ik} + \Gamma^m_{in}\Gamma^n_{jk} - \Gamma^m_{jn}\Gamma^n_{ik})g_{lm}$.

A Riemannian manifold $M$ is said to have negative curvature if $K(X,Y) \leq 0$ for all linearly independent tangent vectors $X,Y \in T_p$ at every point $p \in M$ and have strictly negative curvature if $K(X,Y) \leq -\delta$ ($\delta$ is an arbitrary positive number) for all linearly independent tangent vectors $X,Y \in T_p$ at every point $p \in M$.

Thus, we have the following theorem.

THEOREM 3. *Let $(M,g)$ be a complete simply connected model manifold endowed with the Fisher metric. If $(M,g)$ has strictly negative curvature ($d = 2$) or has negative curvature ($d \geq 3$), then there exist Bayesian predictive distributions asymptotically dominating the Bayesian predictive distribution based on the Jeffreys prior.*

Note that some global differential geometric properties of the model manifold are essential in the present theory, although the field of information geometry has hitherto required only the theory of the locally characterizable properties of manifolds; see [3], page 1.

**3. Conformal transformations corresponding to prior distributions.** We investigate constructing methods for shrinkage priors asymptotically dominating various kinds of vague priors other than the Jeffreys prior. For instance, right invariant priors are more recommended for group models than the Jeffreys priors; see [20, 29]. In this section we assume that the dimension $d$ of the model manifold $M$ is greater than 2.

We introduce conformal transformations corresponding to prior densities and show that there exist shrinkage predictive distributions if the model manifold endowed with the conformally transformed metric satisfies some differential geometric conditions.

A transformation of the metric tensor $g_{ij}(\theta)$ of the form $\tilde{g}_{ij}(\theta) = \nu(\theta) g_{ij}(\theta)$, where $\nu(\theta)$ is a positive function on $M$, is called a conformal transformation. Refer to [11] for details of conformal transformations.

From (3), the difference between the risks of the Bayesian predictive distributions based on prior densities $f$ and $h$ is

$$N^2\{E[D(p(y|\theta), p_h(y|x^{(N)}))|\theta] - E[D(p(y|\theta), p_f(y|x^{(N)}))|\theta]\}$$
$$= \left(\frac{h}{\pi_J}\right)^{2/(d-2)} \left[\frac{1}{2}\left(\frac{h}{f}\right)^2 \tilde{g}^{ij} \partial_i \frac{f}{h} \partial_j \frac{f}{h} - \frac{h}{f}\tilde{\Delta}\frac{f}{h}\right] + o(1),$$

where $\tilde{\Delta}$ denotes the Laplacian corresponding to the metric $\tilde{g}_{ij}(\theta) := \{h(\theta)/\pi_J(\theta)\}^{2/(d-2)} g_{ij}(\theta)$. Thus, we obtain the following theorem in the same way as in the proof of Theorem 1.

THEOREM 4. *Let $h(\theta)$ be a smooth prior density and let $(M, \tilde{g})$ be the model manifold endowed with the metric defined by $\tilde{g}_{ij}(\theta) := \{h(\theta)/\pi_J(\theta)\}^{2/(d-2)} g_{ij}(\theta)$, where $\pi_J(\theta)$ is the density of the Jeffreys prior. The dimension of $M$ is assumed to be greater than 2.*

*If there exists a Green function $\tilde{G}(\xi, \theta)$ on the Riemannian manifold $(M, \tilde{g})$, there exist Bayesian predictive distributions asymptotically dominating the Bayesian predictive distribution based on the prior density $h(\theta)$. In particular, the Bayesian predictive distribution based on the prior density $\pi_{\tilde{G}}(\theta) d\theta := \tilde{G}(\xi, \theta) h(\theta) d\theta$, where $\xi \in M$ is an arbitrary fixed point, asymptotically dominates the Bayesian predictive distribution based on the prior density $h(\theta)$.*

In the same way we proved Theorem 2 and 3, we have the following theorems.

THEOREM 5. *Let $f(\theta)$ and $h(\theta)$ be smooth prior densities on a model manifold $M$ ($d \geq 3$). The Bayesian predictive distribution based on $f(\theta)$ asymptotically dominates the Bayesian predictive distribution based on $h(\theta)$ if and only if $(f/h)^{1/2}$ is a nonconstant positive superharmonic function on* the model manifold $(M, \tilde{g})$ endowed with the conformally transformed metric $\tilde{g}_{ij}(\theta) := \{h(\theta)/\pi_{\mathrm{J}}(\theta)\}^{2/(d-2)} g_{ij}(\theta)$.

THEOREM 6. *If a model manifold $(M, \tilde{g})$ $(d \geq 3)$ endowed with the conformally transformed metric $\tilde{g}_{ij}(\theta) := \{h(\theta)/\pi_{\mathrm{J}}(\theta)\}^{2/(d-2)} g_{ij}(\theta)$ is complete, simply connected and has negative curvature, there exist Bayesian predictive distributions asymptotically dominating the Bayesian predictive distribution based on the prior density $h(\theta)$.*

**4. Examples.** In this section we see several examples of Bayesian predictive distributions based on shrinkage priors constructed by using the methods introduced in the previous sections.

In estimation theory, it is known that asymptotic domination of one estimator over another does not always means exact finite-sample domination because there are examples where the convergence of the risk expansion is not uniform over the parameter space (see [14, 15, 23, 24, 25]), and further studies that bridge asymptotic and exact theories are required. The same difficulty exists in asymptotic prediction theory.

Nevertheless, in the following examples, many Bayesian predictive distributions based on shrinkage priors constructed by using the asymptotic theoretical methods exactly dominate Bayesian predictive distributions based on vague priors. Therefore, the methods introduced in the previous sections are useful tools to construct shrinkage predictive distributions for practical use.

EXAMPLE 1 (*The multivariate normal model with a known covariance matrix*). We consider the $d$-dimensional Normal model $\mathrm{N}_d(\mu, \Sigma)$, where $\mu = (\mu_1, \mu_2, \ldots, \mu_d) \in \mathbb{R}^d$ is an unknown mean vector and $\Sigma$ is a known variance-covariance matrix. We consider the problem of predicting $y \sim \mathrm{N}_d(\mu, \Sigma)$ using $x^{(N)}$, that is, a set of $N$ independent observations from the same density.

The $(i,j)$-component of the Fisher information matrix $g_{ij}$ does not depend on $\mu$. We assume that $g_{ij} = 1$ for $i = j$ and $g_{ij} = 0$ for $i \neq j$ without loss of generality. Thus, the model manifold $(M, g)$ endowed with the Fisher metric is isometric to $d$-dimensional Euclidean space.

The Jeffreys prior $\pi_{\mathrm{J}}(\mu) \propto 1$, which is invariant under the translation group, is commonly used as a vague prior for $\mu$. The Bayesian predictive density $p_{\pi_{\mathrm{J}}}(y|x^{(N)})$ based on $\pi_{\mathrm{J}}(\mu)$ is the best predictive density that is invariant under the translation group.

The minimal positive fundamental solution of the heat equation $\partial u(t, \theta)/\partial t = \Delta u(t, \theta)$ on $\mathbb{R}^d$ endowed with the usual Euclidean metric is given by $\gamma_t(\xi, \mu) = (4\pi t)^{-d/2} \exp\{-\|\mu - \xi\|^2/(4t)\}$, where $\mu, \xi \in \mathbb{R}^d$.

When $d \leq 2$, the integral (1) becomes infinite and shrinkage priors do not exist. This fact corresponds to the relation discussed by Brown [6] between the recurrence properties of Brownian motion on $\mathbb{R}^d$ and the existence of shrinkage estimators for the multivariate normal model with known covariance matrix. When $d \geq 3$, the integral, which is the Green function on $d$-dimensional Euclidean space, is given by $G(\xi, \mu) = \{\Gamma(d/2-1)/(4\pi^{d/2})\}\|\mu - \xi\|^{-(d-2)}$.

The Green prior defined by $\pi_G(\mu)\,d\mu = G(0,\mu)\pi_J(\mu)\,d\mu \propto \|\mu\|^{-(d-2)}\,d\mu$ coincides with Stein's prior $\pi_S(\mu)$ [26]. By Theorem 1, the Bayesian predictive distribution $p_G(y|\bar{x})$ based on Stein's prior asymptotically dominates $p_{\pi_J}(y|\bar{x})$. These asymptotic results also hold for general multivariate location models.

The explicit form $p_{\pi_G}(y|x^{(N)})$ for the $d$-dimensional Normal model was obtained and it was shown that $p_{\pi_G}(y|x^{(N)})$ exactly dominates $p_{\pi_J}(y|x^{(N)})$ for arbitrary $N > 0$; see [19]. Recently, George, Liang and Xu [13] showed the Bayesian predictive distribution based on a prior density $f$ exactly dominates $p_{\pi_J}(y|x^{(N)})$ when $N$ is sufficiently large if and only if $\sqrt{f}$ is a positive superharmonic function. This result for the multivariate Normal model corresponds to Theorem 2.

Next, we consider the conformal transformation corresponding to the Green prior $\pi_G(\mu)\,d\mu$. Here we assume that the parameter space is $\Theta = \mathbb{R}^d - \{O\}$ for simplicity. Then the model manifold $M$ is homeomorphic to the "cylinder" $S^{d-1} \times \mathbb{R}$, where $S^{d-1}$ denotes the $(d-1)$-dimensional unit sphere. The conformal transformation corresponding to the prior $\pi_G$ is given by $\tilde{g}_{ij} = (\pi_G(\mu)/\pi_J(\mu))^{2/(d-2)} g_{ij} = (\sum \mu_i^2)^{-1} g_{ij}$. The Riemannian manifold $(M, \tilde{g})$ can be imbedded in Euclidean space $\mathbb{R}^{d+1}$ endowed with the usual metric by the map $(\mu_1, \mu_2, \ldots, \mu_d) \mapsto ((\sum_{i=1}^d \mu_i^2)^{-1/2}\mu_1, \ldots, (\sum_{i=1}^d \mu_i^2)^{-1/2}\mu_d, (1/2) \times \log(\sum_{i=1}^d \mu_i^2))$.

There does not exist a Green function on the Riemannian manifold $(M, \tilde{g})$, because the integral (1) becomes infinite. Thus, a predictive distribution asymptotically dominating $p_{\pi_G}(y|x^{(N)})$ based on the Green prior cannot be constructed by using the method introduced in Section 3. This fact seems to be related to the admissibility of the shrinkage predictive distribution.

EXAMPLE 2 (*Location-scale models*). Let $p(x)$ be a probability density on $\mathbb{R}$ that is symmetric about the origin. We consider the location-scale model $p(x|\mu, \sigma)\,dx := (1/\sigma)p((x - \mu)/\sigma)\,dx$, where $\mu \in \mathbb{R}$ and $\sigma > 0$ are unknown parameters. Without loss of generality, we can assume that the metric tensor coefficients are given by $g_{\mu\mu} = a/\sigma^2$, $g_{\sigma\sigma} = a/\sigma^2$ and $g_{\mu\sigma} = 0$ by rescaling $\mu$, where $a > 0$ is a constant depending on $p(x)$. The model manifold is the hyperbolic plane $H^2(-1/a)$ with constant curvature $K = -1/a$. The Laplacian on the model manifold is given by $\Delta = (\sigma^2/a)(\partial^2/\partial\mu^2 + \partial^2/\partial\sigma^2)$. The Green function on $H^2(-1/a)$ is given by $G((\mu, \sigma), (0, 1)) = -1/(2\pi) \log \tanh(\rho(\mu, \sigma)/(2\sqrt{a}))$, where $\rho(\mu, \sigma) = \text{dis}((0, 1), (\mu, \sigma))$; see, for example, [4, 9].

Thus, the Bayesian predictive distribution based on the prior $\pi_{\mathrm{G}}(\mu, \sigma) \, d\mu \, d\sigma \propto G((\mu, \sigma), (0, 1)) \pi_{\mathrm{J}}(\mu, \sigma) \, d\mu \, d\sigma$ asymptotically dominates the Bayesian predictive distribution based on the Bayesian predictive distribution based on the Jeffreys prior.

The location-scale model is a group model. The Jeffreys prior $\pi_{\mathrm{J}}(\mu, \sigma) \propto 1/\sigma^2$ is the left invariant prior. However, the best invariant predictive distribution is the Bayesian predictive distribution based on the right invariant prior $\pi_{\mathrm{R}}(\mu, \sigma) \propto 1/\sigma$; see [20, 29]. Here $\pi_{\mathrm{R}}/\pi_{\mathrm{J}} \propto \sigma$ is a positive harmonic function on the model manifold and satisfies the condition of Theorem 2. Furthermore, there exist Bayesian predictive distributions based on positive superharmonic priors asymptotically dominating the best invariant predictive distribution. The details of this topic will be discussed in another paper.

EXAMPLE 3 (*The $2 \times 2$ Wishart model $W_2(m, \Sigma)$*). Suppose that we have a set of independent observations $X(1), X(2), \ldots, X(N)$ from the $2 \times 2$ Wishart distribution $W_2(m, \Sigma)$ with $m$ degrees of freedom. The density of the $2 \times 2$ Wishart distribution $W_2(m, \Sigma)$ is given by

$$p(X|\Sigma) \, dX = \frac{1}{2^m \Gamma_2(m/2) |\Sigma|^{m/2}} |X|^{(m-3)/2} \exp\left(-\frac{1}{2} \mathrm{tr} \Sigma^{-1} X\right) dX,$$
$$X > 0, m \geq 2.$$

Then the distribution of the sufficient statistic $\tilde{X} := \sum_{l=1}^{N} X(l)$ is the $2 \times 2$ Wishart distribution $W_2(Nm, \Sigma)$ with $Nm$ degrees of freedom.

Let $(M, g)$ be a Riemannian manifold composed of $2 \times 2$ positive definite matrices endowed with the Fisher metric. The inner product of tangent vectors $A$ and $B$ at a point $\Sigma \in M$ is given by $\frac{m}{2} \mathrm{tr}(\Sigma^{-1} A \Sigma^{-1} B)$, and the Jeffreys prior is given by $\pi_{\mathrm{J}}(\Sigma) \, d\Sigma \propto |\Sigma|^{-3/2} \, d\Sigma = |\Sigma|^{3/2} \, d\Sigma^{-1}$.

The posterior distribution with respect to the Jeffreys prior is the inverted Wishart distribution $W_d^{-1}(Nm, X)$. The Bayesian predictive distribution for $Y \sim W_2(m, \Sigma)$ based on the Jeffreys prior is given by

$$p_{\pi_{\mathrm{J}}}(Y|\tilde{X}) \, dY = \frac{\Gamma((Nm+m)/2) \Gamma((Nm+m-1)/2)}{\pi^{1/2} \Gamma(Nm/2) \Gamma((Nm-1)/2) \Gamma(m/2) \Gamma((m-1)/2)}$$
$$\times \frac{|\tilde{X}|^{Nm/2} |Y|^{(m-3)/2}}{|\tilde{X} + Y|^{(Nm+m)/2}} \, dY.$$

We construct a shrinkage predictive distribution asymptotically dominating $p_{\pi_{\mathrm{J}}}(Y|\tilde{X})$.

We parameterize $\Sigma$ by

$$(6) \qquad \Sigma = e^\lambda \begin{pmatrix} \cos\frac{\theta}{2} & -\sin\frac{\theta}{2} \\ \sin\frac{\theta}{2} & \cos\frac{\theta}{2} \end{pmatrix} \begin{pmatrix} e^\rho & 0 \\ 0 & e^{-\rho} \end{pmatrix} \begin{pmatrix} \cos\frac{\theta}{2} & \sin\frac{\theta}{2} \\ -\sin\frac{\theta}{2} & \cos\frac{\theta}{2} \end{pmatrix},$$

where $\lambda \in \mathbb{R}$, $\rho > 0$, $0 \le \theta < 2\pi$. The Fisher metric is $g = m(d\lambda^2 + d\rho^2 + \sinh^2\rho\, d\theta^2)$, and the density of the Jeffreys prior is $\pi_{\mathrm{J}}(\lambda, \rho, \theta)\, d\lambda\, d\rho\, d\theta = \pi_{\mathrm{J}}(\Sigma)|(\partial\Sigma/\partial(\rho, \theta, \lambda))|\, d\lambda\, d\rho\, d\theta \propto \sinh\rho\, d\lambda\, d\rho\, d\theta$.

Let $S_{\lambda_0}$ be a submanifold of $(M, g)$ specified by $\lambda = \lambda_0$, where $\lambda_0$ is a constant. The induced metric on the submanifold $S_{\lambda_0}$ is $g = m(d\rho^2 + \sinh^2\rho\, d\theta^2)$. Thus, the Riemannian submanifold $(S_{\lambda_0}, g)$ of the model manifold $(M, g)$ is isometric to the hyperbolic plane $H^2(-1/m)$ with constant curvature $K = -1/m$. Geometric properties of the hyperbolic plane are widely known; see, for example, [9]. The Laplacian on $(M, g)$ is

$$\Delta = \frac{1}{m}\left(\frac{\partial^2}{\partial \lambda^2} + \frac{\partial^2}{\partial \rho^2} + \frac{\cosh\rho}{\sinh\rho}\frac{\partial}{\partial\rho} + \frac{1}{\sinh^2\rho}\frac{\partial^2}{\partial\theta^2}\right).$$

The Riemannian geometric structure of $S_{\lambda_0}$ does not depend on the value of $\lambda_0$. We identify $(S_\lambda, g)$ with $H^2(-1/m)$. The set of submanifolds $\{S_\lambda | \lambda \in \mathbb{R}\}$ is a foliation of the model manifold $(M, g)$.

By Theorem 3, a shrinkage prior dominating the Jeffreys prior on the model manifold $(M, g)$ exists since $(M, g)$ has negative curvature and the dimension $d$ is 3. Here we introduce a shrinkage prior based on the Green function on $(S_\lambda, g)$, which is different from the Green prior $\pi_{\mathrm{G}}$ based on the Green function on $(M, g)$. The Green function on $(S_\lambda, g)$ is given by $\bar{G}_\lambda(e^\lambda I, \Sigma) = -(1/2\pi)\log\{\tanh(\rho/2)\}$, where $|\Sigma| = \exp(2\lambda)$.

We define a function on $(M, g)$ by $h(\Sigma) := \bar{G}_{(1/2)\log|\Sigma|}(|\Sigma|^{1/2}I, \Sigma)$. The function $h(\Sigma)$ is superharmonic on $(M, g)$ and satisfies $\Delta h(\Sigma) = 0$ if $\rho \ne 0$.

We introduce a shrinkage prior distribution defined by $\pi_{\mathrm{S}}(\Sigma)\, d\rho\, d\theta\, d\lambda \propto h(\Sigma)\pi_{\mathrm{J}}(\lambda, \rho, \sigma)\, d\lambda\, d\rho\, d\theta \propto -(1/2\pi)\log\{\tanh(\rho/2)\}\sinh\rho\, d\lambda\, d\rho\, d\theta$. This prior "shrinks" the posterior to the submanifold of $(M, g)$ specified by $\rho = 0$.

By a discussion similar to the proof of Theorem 1, we can see that the Bayesian predictive distribution $p_{\pi_{\mathrm{S}}}(Y|\tilde{X})$ based on $\pi_{\mathrm{S}}$ asymptotically dominates $p_{\pi_{\mathrm{J}}}(Y|\tilde{X})$.

In fact, the explicit Bayesian predictive distribution $p_{\pi_{\mathrm{S}}}(Y|\tilde{X})$ with respect to the prior $\pi_{\mathrm{S}}$ exactly dominates $p_{\pi_{\mathrm{J}}}(Y|\tilde{X})$ based on the Jeffreys prior. The details and some other priors will be discussed in a forthcoming paper.

**Acknowledgments.** The author is grateful to Professor Shinto Eguchi and Professor Tatsuya Kubokawa for helpful comments and to the referee for helpful suggestions.

DEPARTMENT OF MATHEMATICAL INFORMATICS
GRADUATE SCHOOL OF INFORMATION SCIENCE AND TECHNOLOGY
UNIVERSITY OF TOKYO
7-3-1 HONGO, BUNKYO-KU, TOKYO 113-8656
JAPAN
E-MAIL: komaki@mist.i.u-tokyo.ac.jp